\documentstyle[12pt]{article}
\date{ }
\newcommand{\vs}{\vspace{.20in}}

\newcommand{\noin}{\noindent}
\newcommand{\n}{\indent}

\begin{document}


\begin{center} {\bf Two Interesting Properties of the Exponential Distribution} \vspace{.10in}

Robert W. Chen \vspace{.20in}

{\bf Abstract} \end{center} \vspace{.05in}

Let $X_1, X_2,\ldots, X_n$ be $n$ independent and identically distributed random variables, here
$n \geq 2.$ Let $X_{(1)}, X_{(2)}, \ldots, X_{(n)}$ be the order statistics of  $X_1, X_2,..., X_n.$
In this note we proved that: (I) If $X_1, X_2,..., X_n$ are exponential random variables with parameter
$c > 0,$ then the "correlation coefficient" between $X_{(k)}$ and $X_{(k+t)}$ is strictly increasing in
$k$ from $1$ to $m,$ and then is strictly decreasing in $k$ from $m$ to $n - t,$ here $t$ is a fixed
integer between $1$ and $n - 3,$ and $m = (n - t)/2$ if $n - t$ is even, $m = (n - t + 1)/2$ if $n - t$
is odd. We also proved that if $t = n - 2$, then the "correlation coefficient" between $X_{(1)}$ and
$X_{(n-1)}$ is greater than the "correlation coefficient" between $X_{(2)}$ and$X_{(n)}.$ (II) The
"correlation coefficient" between $X_{(k)}$ and $X_{(k+t)}$ for the exponential random variables is always
less than the "correlation coefficient" between $X_{(k)}$ and $X_{(k+t)}$ for the uniform random variables
for all $k$ and $t$ such that $k + t \leq n.$ A combinatorial identity is also given as a bi-product. \vs

\noin MATHEMATICS SUBJECT CLASSIFICATION (2000): Primary 62E99,\\ Secondary 62F99. \\
\\~
Key words and phases: Exponential distribution, order statistics, correlation \\
coefficient, combinatorial identity.\\
\\~
\noin Robert W. Chen: Department of Mathematics, University of Miami, \\
\noin1365 Memorial Drive, Coral Gables, FL 33124-4250. \\
\\~
email: chen@math.miami.edu \\

\newpage

Let $X_1, X_2,..., X_n$ be $n$ independent and identically distributed exponential random variables
with parameter $c > 0,$ here $n \geq 2.$ Let $X_{(1)}, X_{(2)},.., X_{(n)}$ be the order statistics
of $X_1, X_2,..., X_n.$ Without loss of generality, we can and assume that $c = 1.$ The joint
probability density function of $X_{(1)}, X_{(2)},..., X_{(n)}$ is $f(x_{(1)}, x_{(2)},..., x_{(n)})
= n!exp\{-[x_{(1)}+x_{(2)}+...+x_{(n)}]\},$ here $0 \leq x_{(1)} \leq x_{(2)} \leq ... \leq x_{(n)} <
\infty.$ Now let $X_{(1)} = Y_1,$ $X_{(2)} = Y_1 + Y_2,...,$ $X_{(n)} = Y_1 + Y_2 + ... +Y_n.$ Then the
joint probability density function of $Y_1, Y_2,..., Y_n$ is $g(y_1,y_2,..,y_n) =
ne^{-ny_{1}}(n-1)e^{-(n-1)y_{2}}..2e^{-2y_{n-1}}e^{-y_{n}},$ where $0 \leq y_i < \infty$ for all $i =1,
2, .., n.$ It is easy to see that $Y_1, Y_2,..., Y_n$ are mutually independent and $Y_i$ is an exponential
random variable with parameter $1/(n+1-i)$ for all $i =1, 2, ..., n.$ Since $X_{(k)} = \sum _{i=1}^k
Y_i,$ $E(X_{(k)}) = E(\sum _{i=1}^kY_i) = \sum_{i=1}^k\frac {1}{(n+1-i)} = \sum_{i=n+1-k}^n\frac
{1}{i}$ and $Var(X_{(k)}) = Var(\sum _{i=1}^kY_i) = \sum_{i=1}^k\frac {1}{(n+1-i)^2} = \sum_{i=n+1-k}^n\frac
{1}{i^2}$ for all $k = 1, 2, ..., n.$ Also $Cov(X_{(k)}, X_{(k+t)}) = Cov(\sum _{i=1}^kY_i, \sum _{i=1}^kY_i
+ \sum _{i=k+1}^{k+t}Y_i)$ $= Cov(\sum _{i=1}^kY_i, \sum _{i=1}^kY_i) + Cov(\sum _{i=1}^kY_i, \sum _{i=k+1}^{k+t}Y_i)$ $= Var(X_{(k)}) + 0 = Var(X_{(k)})$ since $\sum _{i=1}^kY_i$ and $\sum _{i=k+1}^{k+t}Y_i$
are independent, here $1 \leq t \leq n - k.$ \vs

$E(X_{(n)}^2) = \int_0^\infty x^2n(1-e^{-x})^{n-1}e^{-x}dx $ $= n\sum_{j=0}^{n-1}{n-1 \choose j}
(-1)^j\frac {2}{(j+1)^3} = 2[\sum_{i=1}^n{n \choose i}$ $\frac {(-1)^{i+1}}{i^2}]$ $ = Var(X_{(n)}) +
[E(X_{(n)})]^2$ $= \sum_{i=1}^n\frac {1}{i^2} + [\sum_{i=1}^n\frac {1}{i}]^2.$ Therefore, we have
the following combinatorial identity

$$ (1)\;\;\sum_{i=1}^n{n \choose i}\frac {(-1)^{i+1}}{i^2} = \frac {1}{2}\{\sum_{i=1}^n\frac {1}
{i^2} + [\sum_{i=1}^n\frac {1}{i}]^2\}.$$

\noin The following combinatorial identity

$$(2)\;\;\sum_{i=1}^n{n \choose i}(-1)^{i+1}\frac {1}{i} = \sum_{i=1}^n\frac {1}{i}.$$

\noin is known. However, it can be derived simply by computing $E(X_{(n)}) = \int _0^\infty
nx(1-e^{-x})^{n-1}e^{-x}dx $ $= n\sum_{j=0}^{n-1}{n-1 \choose j}(-1)^j\frac {1}{(j+1)^2}$
$= \sum_{i=1}^n{n \choose i}(-1)^{i+1}\frac {1}{i} = \sum_{i=1}^n\frac {1}{i}$ since $E(X_{(n)}) = \sum_{i=1}^n\frac {1}{i}.$ \vs

\noin The combinatorial identity (1) might be new. \vs

Let $\rho _{k,t}$ be the "correlation coefficient" between $X_{(k)}$ and $X_{(k+t)}$, where $1 \leq k
\leq n - t$ and $t$ is a fixed positive integer such that $1 \leq t$ and $k + t \leq n.$ First we will
prove that $\rho _{k,t}$ is strictly increasing in $k$ from $1$ to $m$ and then is strictly decreasing
in $k$ from $m$ to $n - t.$ It is easy to check that $\rho^2_{k,t} = [\sum_{n+1-k}^n \frac{1}{i^2}]/[\sum_{n+1-k-t}^n\frac {1}{i^2}]$ since $Cov(X_{(k)}, X_{(k+t)}) = Var(X_{(k)}).$ Since $t$
is fixed, we will let $h(k) = \rho^2_{k,t}$ and are interested in the function $h(k)$ for $k$ from $1$ to $n - t.$  First we give a few examples. \vs

\newpage

\noin Example 1: $\;\; n = 5.$ \\
\\~
$t = 1,\; h(1) \approx 0.390, h(2) \approx 0.480, h(3) \approx 0.461, h(4) \approx 0.317.$ \\
\\~
$t = 2,\; h(1) \approx 0.187, h(2) \approx 0.221, h(3) \approx 0.146.$ \\
\\~
$t = 3,\; h(1) \approx 0.086, h(2) \approx 0.070.$ \\
\\~
$t = 4,\; h(1) \approx 0.027.$ \vs

\noin Example 2: $\;\; n = 6.$ \\
\\~
$t = 1,\; h(1) \approx 0.410, h(2) \approx 0.520, h(3) \approx 0.540, h(4) \approx 0.491, h(5) \approx 0.329.$ \\
\\~
$t = 2,\; h(1) \approx 0.213, h(2) \approx 0.281, h(3) \approx 0.265, h(4) \approx 0.162.$ \\
\\~
$t = 3,\; h(1) \approx 0.115, h(2) \approx 0.138, h(3) \approx 0.087.$ \\
\\~
$t = 4,\; h(1) \approx 0.056, h(2) \approx 0.045.$ \\
\\~
$t = 5, \; h(1) \approx 0.019.$ \vs

\noin Example 3: $\;\; n = 8.$\\
\\~
$t = 1,\; h(1) \approx 0.434, h(2) \approx 0.565, h(3) \approx 0.615, h(4) \approx 0.624, h(5) \approx 0.599,$\\
\\~
\n \n $ h(6) \approx 0.526, h(7) \approx 0.345.$ \\
\\~
$t = 2,\; h(1) \approx 0.245, h(2) \approx 0.347, h(3) \approx 0.384, h(4) \approx 0.374, h(5) \approx 0.315,$\\
\\~
\n \n $ h(6) \approx 0.182.$ \\
\\~
$t = 3,\; h(1) \approx 0.151, h(2) \approx 0.217, h(3) \approx 0.230, h(4) \approx 0.197, h(5) \approx 0.109.$ \\
\\~
$t = 4,\; h(1) \approx 0.094, h(2) \approx 0.130, h(3) \approx 0.121, h(4) \approx 0.068.$ \\
\\~
$t = 5,\; h(1) \approx 0.056, h(2) \approx 0.068, h(3) \approx 0.042.$ \\
\\~
$t = 6,\; h(1) \approx 0.030, h(2) \approx 0.024.$ \\
\\~
$t = 7,\; h(1) \approx 0.010.$ \vs

\noin Example 4: $\;\; n = 9.$\\
\\~
$t = 1,\; h(1) \approx 0.441, h(2) \approx 0.578, h(3) \approx 0.635, h(4) \approx 0.656, h(5) \approx 0.650,$\\
\\~
\n \n$ h(6) \approx 0.617, h(7) \approx 0.537, h(8) \approx 0.351.$ \\
\\~
$t = 2,\; h(1) \approx 0.355, h(2) \approx 0.367, h(3) \approx 0.417, h(4) \approx 0.426, h(5) \approx 0.401,$\\
\\~
\n \n $h(6) \approx 0.331, h(7) \approx 0.188.$ \\
\\~
$t = 3,\; h(1) \approx 0.162, h(2) \approx 0.241, h(3) \approx 0.271, h(4) \approx 0.263, h(5) \approx 0.215,$\\
\\~
\n \n $h(6) \approx 0.116.$ \\
\\~
$t = 4,\; h(1) \approx 0.106, h(2) \approx 0.157, h(3) \approx 0.167, h(4) \approx 0.141, h(5) \approx 0.075.$ \\
\\~
$t = 5,\; h(1) \approx 0.069, h(2) \approx 0.097, h(3) \approx 0.090, h(4) \approx 0.050.$ \\
\\~
$t = 6,\; h(1) \approx 0.043, h(2) \approx 0.052, h(3) \approx 0.031.$ \\
\\~
$t = 7,\; h(1) \approx 0.022, h(2) \approx 0.018.$ \\
\\~
$t = 8,\; h(1) \approx 0.008.$ \vs

From these examples, we can see that for a fixed $"t"$, $h(k)$ is strictly increasing and then strictly decreasing for $k$ from $1$ to $n - t,$ except that when $t = n - 2,$ then $h(1) > h(2)$ (when $t = n - 2,$ $k$ can be $1$ or $2$ only).\vs

\noin Theorem 1: \vs

(I) For any fixed $t$ between $1$ and $n - 3,$ $h(k)$ is strictly increasing for $1 \leq k \leq m$ and
is strictly decreasing for $m \leq k \leq n - t,$ where $m = (n - t)/2$ if $n - t$ is even $= (n+1-t)/2$
if $n - t$ is odd. \vs

(II) For $t = n - 2,$ then $h(1) > h(2)$.\vs

Before we prove "Theorem 1", we state a "Lemma" without a proof since it is easy to
check. \vs

\noin Lemma: Assume that a, b, c, and d are positive numbers, \\
\\~
(a) If $a/(a+b) > c/d,$ then $a/(a+b) > (a+c)/(a+b+d) > c/d.$ \\
\\~
(b) If $a/(a+b) < c/d,$ then $a/(a+b) < (a+c)/(a+b+d).$ \vs

Now we start to prove "Theorem 1". \vs

(I) For a fixed $t$ between $1$ and $n - 3,$ $m \geq 2.$ We first will show that

$$h(m) = [\sum_{i=n+1-m}^n\frac {1}{i^2}]/[\sum_{i=n+1-m-t}^n\frac {1}{i^2}] >
\frac {(n-m-t)^2}{(n-m)^2}.$$

If this is proved, then by the "Lemma", $h(m) > h(m+1) > (n-m-t)^2/(n-m)^2$ and
$h(m+1) > (n-m-t-1)^2/(n-m-1)^2$ since $(n-m-t)^2/(n-m)^2$ is strictly decreasing in $m$
for a fixed $"t".$ By this process, we will have $h(m) > h(m+1) > ... > h(n-t),$ i.e.,
$h(k)$ is strictly decreasing in $k$ from $m$ to $(n - t).$ To prove that $h(m) >
(n-m-t)^2/(n-m)^2,$ let $N = \sum_{i=n+1-m}^n\frac {1}{i^2}$ and $D = \sum_{i=n+1-m-t}^n
\frac {1}{i^2}.$ It is easy to check that

$$ N > [ \frac {1}{2(n+1-m)^2} + \int _{n+1-m}^n \frac {1}{x^2}dx + \frac {1}{2n^2}]$$
$$= \frac {n^2 + (n+1-m)^2 + 2n(n+1-m)(m-1)}{2n^2(n+1-m)^2}.$$

\noin Also it is easy to check that
$$D < \int _{n+1-m-t}^n \frac{1}{(x-\frac {1}{2})^2}dx = \frac {2}{2n+1-2m-2t} - \frac {2}{2n+1}$$
$$= \frac{4(m+t)}{(2n+1)(2n+1-2m-2t)}.$$

\noin To prove that $h(m) = N/D > (n-m-t)^2/(n-m)^2,$ it is sufficient to show that

$$\frac {(2n+1)(2n+1-2m-2t)[n^2+(n+1-m)^2+2n(n+1-m)(m-1)]}{8(m+t)n^2(n+1-m)^2}$$
$$ >\frac {(n-m-t)^2}{(n-m)^2}$$

\noin since $h(m) = N/D >$

$$\frac {(2n+1)(2n+1-2m-2t)[n^2+(n+1-m)^2+2n(n+1-m)(m-1)]}{8(m+t)n^2(n+1-m)^2}.$$

\noin Since all numbers involved are positive numbers, it is sufficient to show that

$$(3)\;\;(2n+1)(2n+1-2m-2t)[n^2+(n+1-m)^2+2n(n+1-m)(m-1)](n-m)^2$$
$$-\;8(m+t)n^2(n+1-m)^2(n-m-t)^2 > 0.$$

There are two cases to be considered:\vs

\noin $(a)\;\; n = 2m + t.$ \vs

Then to prove the inequality $(3)$ is equivalent to prove the following inequality

$$(4)\;(4m+2t+1)(2m+1)(m+t)^2[(2m+t)^2+(m+t+1)^2+2(2m+t)(m+t+1)(m-1)]$$
$$-\;8(m+t)m^2(m+t+1)^2(2m+t)^2 > 0.$$

\noin After simplification, we have the following inequality:

$$(5)\; (16t-14)m^5+(48t^2-8t+1)m^4+(53t^3+24t^2+4t+4)m^3$$
$$+\; (24t^4+24t^3+5t^2+10t+1)m^2+(4t^5+6t^4+2t^3+8t^2+2t)m$$
$$+\; t^2(2t+1) > 0$$

\noin since $t \geq 1$ and $m \geq 2.$ Hence $h(m) >(n-m-t)^2/(n-m)^2$ when $n = 2m+t.$ \vs

\noin $(b)\;\; n = 2m + t - 1.$ \vs

Then to prove the inequality $(3)$ is equivalent to prove the following inequality

$$(6)\;(4m+2t-1)(2m-1)(m+t-1)^2[(2m+t-1)^2+(m+t)^2+2(2m+t-1)(m+t)(m-1)]$$
$$-\;8(m+t)^3(2m+t-1)^2(m-1)^2 > 0.$$

\noin After simplification, the left hand side of the inequality $(6)$ is

$$(7)\; (48t-14)m^5+(144t^2-120t+47)m^4+(156t^3-256t^2+128t-60)m^3$$
$$+\; (72t^4-216t^3+155t^2-80t+36)m^2+(12t^5-74t^4+86t^3-42t^2+28t-10)m$$
$$+\; (-8t^5+16t^4-10T63+5t^2-4t+1).$$

\noin We have to re-arrange $(7)$ by using the fact that $m \geq 2$ to show that

$$(8)\; (48t-14)m^5+(144t^2-120t+47)m^4+(156t^3-256t^2+128t-60)m^3$$
$$+\; (72t^4-216t^3+155t^2-80t+36)m^2+(12t^5-74t^4+86t^3-42t^2+28t-10)m$$
$$+\; (-8t^5+16t^4-10T63+5t^2-4t+1) \geq 14(t-1)m^5+26(t-1)^2m^4$$
$$+\; [67t(t-1)^2+18(t-1)]m^3 + 19t^2(t-1)^2m^2$$
$$+\; (8t^5+32t^4+81t^3+230t^2+98t+62)m +(16t^4+5t^2+1) > 0$$

\noin since $t \geq 1$ and $m \geq 2.$ Hence the inequality $(6)$ holds and
$h(m) > \frac {(n-m-t)^2}{(n-m)^2}.$ \vs

By the "Lemma", we can conclude both cases that $h(m) > h(m+1) >
\frac {(n-m-t)^2}{(n-m)^2} > \frac {(n-m-t-1)^2}{(n-m-1)^2}.$ Therefore,
$h(m+1) > h(m+2).$  By this process, we have proved that $h(k)$ is strictly
decreasing in $k$ from $m$ to $n-t,$ here $t$ is a fixed integer between $1$ and
$(n-3).$ \vs

Now we have to prove that
$$h(m-1) = [\sum_{i=n+2-m}^n\frac {1}{i^2}]/[\sum_{i=n+2-m-t}^n\frac {1}{i^2}] <
\frac {(n+1-m-t)^2}{(n+1-m)^2}.$$

\noin As above let $N = \sum_{i=n+2-m}^n\frac {1}{i^2}$ and $D = \sum_{i=n+2-m-t}^n
\frac {1}{i^2}.$ It is easy to see that

$$ N < \int _{n+2-m}^{n+1} \frac {1}{(x-0.5)^2}dx = \frac {2}{(2n+3-2m)} - \frac {2}{(2n+1)}$$
$$= \frac {4(m-1)}{(2n+1)(2n+3-2m)}.$$

\noin Also it is easy to check that
$$D > \frac {[n^2+(n+2-m-t)^2+2n(n+2-m-t)(m+t-2)]}{[2n^2(n+2-m-t)^2]}.$$

\noin To prove that $h(m-1) < \frac {(n+1-m-t)^2}{(n+1-m)^2},$ it is sufficient to show that

$$(9)\; \frac {8(m-1)n^2(n+2-m-t)^2}{(2n+1)(2n+3-2m)[n^2+(n+2-m-t)^2+2n(n+2-m-t)(m+t-2)]}$$
$$<\; \frac {(n+1-m-t)^2}{(n+1-m)^2}.$$

\noin As above, it is equivalent to show that
$$(10)\; (2n+1)(2n+3-2m)[n^2+(n+2-m-t)^2+2n(n+2-m-t)(m+t-2)](n+1-m-t)^2$$
$$-\; 8(m-1)n^2(n+2-m-t)^2(n+1-m)^2 > 0.$$ \vs

There are also two cases to be considered.\vs

\noin $(c)\;\; n = 2m + t.$ \vs

Then to prove the inequality $(10)$ is equivalent to prove the following inequality
$$(11)\;(4m+2t+1)(2m+2t+3)(m+1)^2[(2m+t)^2+(m+2)^2+2(2m+t)(m+2)(m+t-2)]$$
$$-\;8(m-1)(m+2)^2(2m+t)^2(m+t+1)^2 > 0.$$ \vs

\noin After simplification, the inequality $(11)$ becomes
$$(12)\; (48t-14)m^5 + (96t^2+242t-21)m^4 + (60t^3+444t^2+472t+20)m^3$$
$$+\; (12t^4+256t^3+695t^2+290t+59)m^2 + (48t^4+332t^3+282t^2+20t+44)m$$
$$+\; (52t^4+72t^3-t^2+8t+12) > 0$$

\noin since $t \geq 1$ and $m \geq 2.$ \vs

\noin $(d)\;\; n = 2m + t - 1.$ \vs

Then the inequality $(10)$ becomes \vs
$$(13)\; m^2(4m+2t-1)(2m+2t+1)[(m+1)^2+(2m+t-1)^2+2(2m+t-1)(m+1)(m+t-2)]$$
$$-\;8(m-1)(m+2)^2(2m+t)^2(m+t-1)^2 > 0.$$ \vs

\noin After simplification, the inequality $(13)$ becomes
$$(14)\; (16t-14)m^5 + (32t^2+2t+25)m^4 + (20t^3+36t^2+52t-4)m^3$$
$$+\; (4t^4+32t^3+53t^2-56t+2)m^2 + (8t^4+32t^3-56t^2+16t)m + 8t^2(t-1)^2 > 0$$
\noin since $t \geq 1$ and $m \geq 2.$ Hence $h(m-1) < \frac {(n+1-m-t)^2}{(n+1-m)^2}$
and $h(m-1) < h(m).$ \vs

Now we have to show that $h(k)$ is strictly increasing in $k$ from $1$ to $m.$ Suppose
not, then there exists a $k$ such that $h(k)\geq h(k+1),$ where $1 \leq k \leq m-2$ since
$h(m-1) < h(m).$  If $h(k) = h(k+1),$ then $h(k) = \frac {(n-k-t)^2}{(n-k)^2}$ and
$h(k) = h(k+1) > \frac {(n-1-k-t)^2}{(n-1-k)^2}.$ By the "Lemma", then $h(k+1) > h(k+2) >...
> h(m-1) > h(m)$ and we get a contradiction. If $h(k) > h(k+1),$ then $h(k+1) > h(k+2) >...
> h(m-1) > h(m)$ and we get a contradiction again. Hence $h(k)$ is strictly increasing in $k$
from $1$ to $m.$ The part (I) of the "Theorem 1" is proved. \vs

To complete the proof of the "Theorem 1", now we have to prove the part (II) of the "Theorem 1". \vs

When $t = n - 2$ and $n \geq 3,$ $k$ can be $1$ or $2$ only. Now we will show that $h(1) > h(2).$
$h(1) = \frac {1}{n^2\sum_{i=2}^n \frac {1}{i^2}},$ to prove $h(1) > h(2),$ we only need to show
$\frac {1}{n^2\sum_{i=2}^n \frac {1}{i^2}} > \frac {1}{(n-1)^2}.$ It is easy to see that
$$(15) \; \sum_{i=2}^n \frac {1}{i^2} < \int _2^{n+1} \frac {1}{(x-0.5)} = \frac {2}{3} - \frac{2}{2n+1} = \frac {4(n-1)}{3(2n+1)}.$$
\noin If $\frac {4(n-1)}{3(2n+1)} < \frac {(n-1)^2}{n^2},$ then $\sum_{i=2}^n\frac {1}{i^2}
< \frac {(n-1)^2}{n^2}.$ $\frac {4(n-1)}{3(2n+1)} < \frac {(n-1)^2}{n^2}$ if $3(2n+1)(n-1)
- 4n^2 > 0.$ $3(2n+1)(n-1) - 4n^2 = 2(n-1)^2 + (n-3) > 0$ since $n \geq 3.$ Therefore, $h(1)
> h(2)$ and the proof of the "Theorem 1" now is complete. \vs

By the same process, we also get an upper bound for $h(m)$. Hence we have the following
inequality.

$$\frac {(2n+1)(2n+1-2m-2t)[n^2+(n+1-m)^2+2n(n+1-m)(m-1)]}{8n^2(n+1-m)^2(m+t)} < h(m) $$
$$< \frac {8n^2(n+1-m-t)^2m}{(2n+1)(2n+1-2m)[n^2+(n+1-m-t)^2+2n(n+1-m-t)(m+t-1)]}.$$ \vs

Suppose that $t = [nx],$ here $[nx]$ is the largest integer $\leq nx$ and $t \leq n - 3.$ Since $h(m) =
\rho^2_{m,t}$ and by the "Lemma", $h(m) < \frac {(n+1-m-t)^2)}{(n+1-m)^2}$, we have the following
inequality for $\rho_{m,t}$
$$ \frac {(n-m-t)}{(n-m)} < \rho_{m,t} < \frac {(n+1-m-t)}{(n+1-m)}.$$
\noin Substitute $[nx]$ for $t$, if $n - [nx]$ is even, we have the following inequality for $\rho_{m,t}$
$$\frac {(n-[nx])}{(n+[nx])} < \rho_{m,t} < \frac {(n-[nx]+2)}{(n+[nx]+2)}.$$
\noin And if $n - [nx]$ is odd, we have the following inequality for $\rho_{m,t}$
$$\frac {(n-[nx]-1)}{(n+[nx]-1)} < \rho_{m,t} < \frac {(n-[nx]+1)}{(n+[nx]+1)}.$$
\noin Further more if $[nx] = nx$ i.e., $nx$ is an integer, then if $n-nx$ is even we have
we the following inequality for $\rho_{m,t}$
$$\frac {(1-x)}{(1+x)} < \rho_{m,t} < \frac {(1-x+\frac {2}{n})}{(1+x+\frac {2}{n})}.$$
\noin And if $n-nx$ is odd we have the following inequality for $\rho _{m,t}$
$$\frac {(1-x-\frac {1}{n})}{(1+x-\frac {1}{n})} < \rho_{m,t} < \frac {(1-x+\frac {1}{n})}{(1+x+\frac {1}{n})}.$$
\noin If $n$ is large, the lower bound and the upper bound are so close, and $\rho_{m,t} \approx \frac {(1-x)}{(1+x)}.$ \vs

In fact, if we replace $m$ by $k$ for $k$ from $1$ to $n-3$, we have upper and lower bounds for $h(k)$ as
follows:
$$(16)\;\frac {(2n+1)(2n+1-2k-2t)[n^2+(n+1-k)^2+2n(n+1-k)(k-1)]}{8n^2(n+1-k)^2(k+t)} < h(k) $$
$$< \frac {8n^2(n+1-k-t)^2k}{(2n+1)(2n+1-2k)[n^2+(n+1-k-t)^2+2n(n+1-k-t)(k+t-1)]}.$$ \vs

It is very easy to compute the lower and upper bounds for $\rho _{k,t}$ for any $k$ and $t$, here
$t$ is fixed and is between $1$ and $n-3$, $k+t \leq n$ even just with a hand calculator. Also the
both bounds are very close to the exact value of $\rho _{k,t}.$ However, even moderate $n$, it needs
some software like Maple or Mathematica to compute $\rho _{k,t}.$ \vs

Now suppose that $X_1, X_2,\ldots, X_n$ are $n$ independent and identically distributed uniform
random variables over the interval $[0, 1]$ here $n \geq 2.$ Let $X_{(1)}, X_{(2)}, \ldots, X_{(n)}$ be
the order statistics of  $X_1, X_2,..., X_n.$ It is well-known that the "correlation coefficient" between $X_{(k)}$ and $X_{(k+t)}$ is equal to $\sqrt {\frac {k(n+1-k-t)}{(k+t)(n+1-k)}}.$ It is easy to see that
it is strictly increasing in $k$ from $1$ to $m$ and then strictly decreasing in $k$ from $m$ to $n-t$ if
$n-t$ is odd. However, if $n-t$ is even, then it is strictly increasing in $k$ from $1$ to $m$ and then
strictly decreasing in $k$ from $m+1$ to $n-t$. For $k = m$ and $k = m + 1$, they are the same. It is
different from the case for the exponential random variables. Further more the "correlation coefficient"
between $X_{(1)}$ and $X_{(n-1)}$ is greater than the "correlation coefficient" between $X_{(2)}$ and
$X_{(n)}$ for the exponential random variables, but the "correlation coefficient" between $X_{(1)}$ and $X_{(n-1)}$ is equal to the "correlation coefficient" between $X_{(2)}$ and $X_{(n)}$ for the uniform
random variables, both of them are equal to $\sqrt {\frac {2}{n(n-1)}}.$ \vs

From our computation, we observed that the "correlation coefficient" between $X_{(k)}$ and $X_{(k+t)}$ for
the exponential random variables is always less than the "correlation coefficient" between $X_{(k)}$ and $X_{(k+t)}$ for the uniform random variables, here $k,\;t$ are positive integers and $k + t \leq n,\; n
\geq 2$. We have the following theorem.\vs

\noin Theorem 2: \vs

\noin The "correlation coefficient" between $X_{(k)}$ and $X_{(k+t)}$ for
the exponential random variables is always less than the "correlation coefficient" between $X_{(k)}$ and $X_{(k+t)}$ for the uniform random variables, here $k,\;t$ are positive integers and $k + t \leq n,\; n
\geq 2$.\\
\\~
\noin It is sufficient to show that $h(k)$ is less than $\frac {k(n+1-k-t)}{(k+t)(n+1-k)}.$ There are three
cases to discuss.\\
\\~
\noin Case I: When $t = 1.$\\
\\~
$$(17)\;\frac {\sum_{n+1-k}^n\frac {1}{i^2}}{\sum_{n-k}^n\frac {1}{i^2}}<\frac {k(n-k)}{(k+1)(n+1-k)}.$$
\\~
It is equivalent to show that
$$(n-k)^2\sum_{n-k}^n\frac {1}{i^2} < \frac {(k+1)(n+1-k)}{n+1}.$$
\\~
\noin After simplification, we have $(2n+1)(2n+1-2k) - 4(n+1)(n-k) = 2k +1 > 0$ since $k \geq 1.$ \\
\\~
\noin Case II: When $k = 1.$\\
\\~
It is equivalent to show that
$$(18)\;\frac {1}{\sum_{n-t}^n \frac {1}{i^2}} < \frac {n(n-t)}{t+1}.$$\\
\\~
Since $\sum_{n-t}^n \frac {1}{i^2} > \frac {n^2+(n-t)^2+2n(n-t)t}{2n^2(n-t)^2},$ it is sufficient to show
that $n^2+(n-t)^2+2n(n-t)t - 2n(n-t)(t+1) > 0.$ After simplification, we have $t^2 > 0$ since $t \geq 1.$ \\
\\~
\noin From now on we will assume that $k, t \geq 2.$ Hence $n \geq 4.$ Recall that
$$h(k) = \frac {N}{D} = \frac {\sum_{n+1-k}^n \frac {1}{i^2}}{\sum_{n+1-k}^n \frac {1}{i^2}}.$$\\
\\~
$$N = \sum_{n+1-k}^n \frac {1}{i^2} < \frac {(2n+1)(2n+3-2k)+4(n+1-k)^2(k-1)}{(2n+1)(2n+3-2k)(n+1-k)^2}.$$\\
\\~
\noin And
$$D = \sum_{n+1-k-t}^n \frac {1}{i^2} > \frac {n^2+(n+1-k-t)^2+2n(n+1-k-t)(k+t-1)}{2n^2(n+1-k-t)^2}.$$\\
\\~
\noin Hence
$$h(k) < $$
$$\frac {2n^2(n+1-k-t)^2[(2n+1)(2n+3-2k)+4(k-1)(n+1-k)^2]}{(2n+1)(2n+3-2k)(n+1-k)^2[n^2+(n+1-k-t)^2+2n(n+1-k-t)(k+t-1)]}.$$\\
\\~
\noin To show $h(k) < \frac {k(n+1-k-t)}{(k+t)(n+1-k)},$ it is sufficient to show that\\
\\~
$$\frac {k(n+1-k-t)}{(k+t)(n+1-k)} > $$
$$\frac {2n^2(n+1-k-t)^2[(2n+1)(2n+3-2k)+4(k-1)(n+1-k)^2]}{(2n+1)(2n+3-2k)(n+1-k)^2[n^2+(n+1-k-t)^2+2n(n+1-k-t)(k+t-1)]}.$$\\
\\~
\noin From now on we will let $n = k + t + x$, where $x$ is a non-negative integer. To show that
$$\frac {k(n+1-k-t)}{(k+t)(n+1-k)} > $$
$$\frac {2n^2(n+1-k-t)^2[(2n+1)(2n+3-2k)+4(k-1)(n+1-k)^2]}{(2n+1)(2n+3-2k)(n+1-k)^2[n^2+(n+1-k-t)^2+2n(n+1-k-t)(k+t-1)]},$$\\
\\~
\noin it is equivalent to show that
$$k(n+1-k-t)(2n+1)(2n+3-2k)(n+1-k)^2[n^2+(n+1-k-t)^2+2n(n+1-k-t)(k+t-1)]$$
$$- 2(k+t)(n+1-k)n^2(n+1-k-t)^2[(2n+1)(2n+3-2k)+4(k-1)(n+1-k)^2] > 0.$$\\
\\~
Replace $n$ by $k+t+x,$ we obtain the following polynomial in $x$

$$(19)\; p(x) = [8k^2 + 8k(t-1) -8t]x^7+[16k^3+(24+56t)k^2+8(5t^2-2t-5)k - 40(t^2+t)]x^6$$

$$+[8k^4+72(t+1)k^3+(144t^2+188t-6)k^2+(80t^3+36t^2-182t-74)k-(80t^3+176t^2+78t)]x^5$$

$$+[8(3t+5)k^4+(120t^2+296t+114)k^3+(176t^3+480t^2+140t-100)k^2$$
$$+(80t^4+144t^3-278t^2-386t-54)k-(80t^4+304t^3+298t^2+74t)]x^4$$

$$+[(24t^2+100t+78)k^4+(88t^3+436t^2+444t+61)k^3+(104t^4+552t^3$$
$$+492t^2-172t-144)k^2+(40t^5+176t^4-130t^3-655t^2-338t+5)k$$
$$-(40t^5+256t^4+434t^3+240t^2-34t)]x^3$$

$$[(8t^3+80t^2+156t+740)k^4+(24t^4+272t^3+560t^2+276t-21)k^3+(24t^5$$
$$+296t^4+576t^3+52t^2-316t-84)k^2+(8t^6+95t^5+58t^4-440t^3-545t^2$$
$$-110t+27)k-(8t^6+104t^5+294t^4+282t^3+90t^2+6t-4)]x^2$$

$$+[(20t^3+90t^2+108t+34)k^4+(60t^4+276t^3+307t^2+44t-33)k^3$$
$$+(60t^5+266t^4+198t^3-168t^2-160t-28)k^2+(20t^6+64t^5-89t^4-312t^3$$
$$-178t^2+9)k-(16t^6+88t^5+140t^4+78t^3+12t^2-8t-8)]x$$

$$+[(12t^3+34t^2+28t+6)k^4+(36t^4+92t^3+53t^2-12t-9)k^3$$
$$(36t^5+74t^4-2t^3-64t^2-24t)k^2+(12t^6+8t^5-51t^4-66t^3-18t^2-1)k$$
$$-(8t^6+24t^5+22t^4+6t^3-4t^2-8t-4)]$$

$$\geq [8k^2 + 8k(t-1) -8t]x^7+[16k^3+(40t^2+96t+8)k - (40t^2+40t)]x^6$$

$$+[8k^4+72(t+1)k^3+(80t^3+324t^2+194t-86)k-(80t^3+176t^2+78t)]x^5$$

$$+[8(3t+5)k^4+(120t^2+296t+114)k^3+(80t^4+496t^3+682t^2-106t-254)k$$
$$-(80t^4+304t^3+298t^2+74t)]x^4$$

$$+[(24t^2+100t+78)k^4+(88t^3+436t^2+444t+61)k^3$$
$$+(40t^5+384t^4+974t^3+329t^2-682t-283)k-(40t^5+256t^4+434t^3+240t^2-34t)]x^3$$

$$+[(8t^3+80t^2+156t+740)k^4+(24t^4+272t^3+560t^2+276t-21)k^3+(8t^6+143t^5$$
$$+650t^4+712t^3t^3-441t^2-742t-141)k-(8t^6+104t^5+294t^4+282t^3+90t^2+6t-4)]x^2$$

$$+[(20t^3+90t^2+108t+34)k^4+(60t^4+276t^3+307t^2+44t-33)k^3+(20t^6+184t^5+443t^4$$
$$+84t^3-514t^2-320t-47)k-(16t^6+88t^5+140t^4+78t^3+12t^2-8t-8)]x$$

$$+[(12t^3+34t^2+28t+6)k^4+(36t^4+92t^3+53t^2-12t-9)k^3+(12t^6+80t^5+97t^4-70t^3$$
$$-146t^2-48t-1)k-(8t^6+24t^5+22t^4+6t^3-4t^2-8t-4)] > 0$$
\noin since $k, t \geq 2$ and $x$ is a non-negative integer. The proof of Theorem 2 now is complete. \vs

Theorem 1 tells us that the "correlation coefficient" between $X{(k)}$ and $X{(k+t)}$ is largest when
$k = \frac {n-t}{2}$ if $n - t$ is even, and $k = \frac {n-t+1}{2}$ if $n - t$ is odd, also the
"correlation coefficient" between $X{(1)}$ and $X{(n-1)}$ is larger than the "correlation coefficient"
between $X{(2)}$ and $X{(n)}$ for the exponential random variables. From our computation, this theorem
does not hold for the random variables with the probability density function $f(x) = 2x$ for $0 \leq x
\leq 1.$ When $n = 3,$ the "correlation coefficient" between $X{(1)}$ and $X{(2)}$ is less than the
"correlation coefficient" between $X{(2)}$ and $X{(3)}.$ Also when $n = 7$ and $t = 1,$ the "correlation coefficient" between $X{(k)}$ and $X{(k+1)}$ is largest when $k = 4 > \frac {n-1}{2}.$ However, when
$n = 6,$ the "correlation coefficient" between $X{(k)}$ and $X{(k+1)}$ is largest when $k = 3 = \frac
{n-t+1}{2}.$ For the random variables with the probability density function $f(x) = 2(1-x)$ for $0 \leq x
\leq 1,$ Theorem 1 seems to hold. So we have the following conjecture:\vs

\noin Conjecture I:\\
\\~
\noin Theorem 1 holds if the probability density function is strictly decreasing. Theorem 1 must be
modified as the "correlation coefficient" between $X{(k)}$ and $X{(k+t)}$ is largest when $k = m =
\frac {n-t+1}{2}$ when $n-t$ is odd, and the "correlation coefficient" between $X{(k)}$ and $X{(k+t)}$ is largest when $k = m+1 = \frac {n-t}{2}+1$ when $n-t$ is even if the probability density function is strictly increasing. We do not have any idea about the case that the probability density is increasing and then
decreasing. \vs

Suppose that $Y_1, Y_2,\ldots, Y_n$ be $n$ independent and identically distributed negative exponential
random variables, here $n \geq 2.$ Let $Y_{(1)}, Y_{(2)}, \ldots, Y_{(n)}$ be the order statistics of
$Y_1, Y_2,..., Y_n.$ Let $Y_i = - X_{n+1-i}$ for all $i = 1,2,..,n$, then $X_1, X_2,\ldots, X_n$ are $n$ independent and identically distributed  exponential random variables, and  $X_{(1)}, X_{(2)}, \ldots, X_{(n)}$ be the order statistics of  $X_1, X_2,..., X_n.$ The "correlation coefficient" between $Y_{(i)}$ and
$Y_{(i+t)}$ is the same as the "correlation coefficient" between $X_{(n+1-i)}$ and $X_{(n+1-i-t)}.$ So the "correlation coefficient" between $Y_{(i)}$ and $Y_{(i+t)}$ is strictly increasing in $i$ from $1$ to $m$, and
is strictly decreasing in $i$ from $m$ to $n-t,$ if $n - t$ is odd and here $m = \frac {n+1-t}{2}.$ And
the "correlation coefficient" between $Y_{(i)}$ and $Y_{(i+t)}$ is strictly increasing in $i$ from $1$ to $m+1$, and is strictly decreasing in $i$ from $m+1$ to $n-t,$ if $n - t$ is even and here $m = \frac {n-t}{2}.$
The \"correlation coefficient" between $X_{(1)}$ and $X_{(n-1)}$ is larger than the "correlation coefficient" between $X_{(2)}$ and $X_{(n)}$. So the "correlation coefficient" between $Y_{(1)}$ and $Y_{(n-1)}$ is less than the "correlation coefficient" between $Y_{(2)}$ and $Y_{(n)}$. \vs

\noin Theorem 2 tells us the "correlation coefficient" between $X_{(k)}$ and $X_{(k+t)}$ of the exponential
random variables is always less than the "correlation coefficient" between $X_{(k)}$ and $X_{(k+t)}$ of the
uniform random variables. From our computation of a few continuous random variables, it looks like this 
property seems to hold. So we make the following conjecture. \vs

\noin Conjecture II:\\
\\~
\noin Theorem 2 holds for any continuous random variables, i.e., the "correlation coefficient" between $X_{(k)}$ and $X_{(k+t)}$ of any continuous random variables is always less than the "correlation coefficient" between $X_{(k)}$ and $X_{(k+t)}$ of the uniform random variables, here $k,\;t$ are positive integers and $k + t \leq n,\; n \geq 2$.\\

\end{document}